\documentclass{article}

\usepackage{enumerate}
\usepackage{amsmath,xspace,amssymb,mathrsfs}

\input xy
\xyoption{all}
\xyoption{2cell}
\UseAllTwocells
\CompileMatrices

\title{Non-canonical isomorphisms}
\author{Stephen Lack%
\thanks{The support of the Australian Research Council and
DETYA is gratefully acknowledged.}
\\School of Computing and Mathematics\\
University of Western Sydney\\
Locked Bag 1797 Penrith South DC NSW 1797\\
Australia\\
email: {\tt s.lack@uws.edu.au}}
\date{}

\newcommand{\A}{{\ensuremath{\mathscr A}}\xspace}
\newcommand{\B}{{\ensuremath{\mathscr B}}\xspace}
\newcommand{\C}{{\ensuremath{\mathscr C}}\xspace}
\newcommand{\D}{{\ensuremath{\mathscr D}}\xspace}
\renewcommand{\phi}{\varphi}

\newcommand{\x}{\times}
\newcommand{\ox}{\otimes}

\newtheorem{theorem}{Theorem}
\newtheorem{lemma}[theorem]{Lemma}
\newtheorem{proposition}[theorem]{Proposition}
\newtheorem{corollary}[theorem]{Corollary}
\newtheorem{preremark}[theorem]{Remark}
\newenvironment{remark}{\begin{preremark}\rm}{\end{preremark}}

\newcommand{\proof}{\noindent{\sc Proof:}\xspace}
\def\endproof{~\hfill$\Box$\vskip 10pt}

\begin{document}

\label{firstpage}
\maketitle

\begin{abstract}
We give two examples of categorical axioms asserting that a
canonically defined natural transformation is invertible where the
invertibility of any natural transformation implies that the canonical
one is invertible. The first example is distributive categories, the
second (semi-)additive ones. We show that each follows from a general 
result about monoidal functors.
\end{abstract}

In any category \D with finite products and coproducts there is a natural
family of maps 
$$\xymatrix @C3pc {
X\x Y+X\x Z \ar[r]^-{\delta_{X,Y,Z}} & X\x(Y+Z) }$$
induced, via the universal property of the coproduct $X\x Y+X\x Z$,
by the morphisms $X\x i$ and $X\x j$, where $i$ and $j$ denote 
the coproduct injections of $Y+Z$. Such a \D is said to be {\em distributive}
\cite{ext, Cockett-dist} if the canonical maps 
are invertible; in other words, if the functor $X\x-:\D\to\D$ preserves 
binary coproducts, for all objects $X$.
As observed by Cockett \cite{Cockett-dist}, it follows that $X\x 0\cong 0$, 
so that $X\x-$ in fact preserves finite coproducts.

Claudio Pisani has asked whether the existence of {\em any} natural family of
isomorphisms
$$\xymatrix @C3pc {
X\x Y+X\x Z \ar[r]^-{\psi_{X,Y,Z}} & X\x(Y+Z)}$$
might imply that \D is distributive. Such $\psi$ are the non-canonical isomorphisms of the title. He suggested that this was probably not the 
case, and this was also my immediate reaction. But in fact it is true!
This is the first result of the paper. 

The second result is an analogue for semi-additive categories. Recall
that a  category is pointed when it has an initial object which is
also  terminal ($1=0$), and that 
for any any two objects $Y$ and $Z$ in a pointed category there is a unique morphism from $Y$ to $Z$ which factorizes through the zero object; this morphism
is called $0_{Y,Z}$ or just $0$. If the category has finite products and 
coproducts, then there is a natural family of morphisms 
$$\xymatrix{
Y+Z \ar[r]^{\alpha_{Y,Z}} & Y\x Z}$$
induced by the identities on $Y$ and $Z$ and the zero morphisms $0:Y\to Z$
and $0:Z\to Y$. The category is semi-additive when these $\alpha_{Y,Z}$ are 
invertible \cite[VII.2]{CWM}. A semi-additive category admits a
canonical enrichment over commutative monoids; conversely, any
category enriched over commutative monoids which has either finite
products or finite coproducts is semi-additive. Our result for 
semi-additive categories asserts once again that the existence of any
natural isomorphism $Y+Z\cong Y\x Z$ implies that the category is
semi-additive.

We also show that the common part of the two arguments follows from a 
general result about monoidal functors; since the individual results
are so easy to prove, however, we give them first, in
Sections~\ref{sect:distributive}
and~\ref{sect:additive} respectively, before turning to the general
result
in Section~\ref{sect:monoidal}.

\section{Non-canonical distributivity isomorphisms}\label{sect:distributive}

This section involves, as in the introduction, a category \D with
finite products
and coproducts and a natural family of isomorphisms $X\x Y+X\x Z\cong
X\x(Y+Z)$. First we show, in the following Lemma, that such a \D will
be distributive if $X\x 0\cong 0$. Later on, we shall see that this
Lemma follows from a more
general  result about coproduct-preserving functors due to Caccamo and
Winskel; and that this in turn is a special case of a still more
general result about monoidal functors: this is our Theorem~\ref{thm:monoidal} below.

\begin{lemma}\label{lemma:pres+}
Suppose that as above that we have natural isomorphisms 
$$\xymatrix{
 X\x Y+X\x Z\ar[r]^-{\psi} & X\x(Y+Z)}$$
and that $X\x 0\cong 0$. Then the category \D is distributive.
\end{lemma}

\proof
If $X\x 0\cong 0$, then $\phi_{X,Y,0}$ gives an isomorphism $X\x Y+X\x 0\cong X\x(Y+0)$, which we can regard as simply being an isomorphism $X\x Y\cong X\x Y$.
By naturality, the diagram
$$\xymatrix{
X\x Y \ar[r]^{\psi_{X,Y,0}} \ar[d]_{i} & X\x Y \ar[d]^{X\x i} \\
X\x Y+X\x Z \ar[r]_-{\psi_{X,Y,Z}} & X\x(Y+Z) }$$ 
commutes, and similarly we have a commutative diagram
$$\xymatrix{
X\x Z \ar[r]^{\psi_{X,0,Z}} \ar[d]_{i} & X\x Y \ar[d]^{X\x i} \\
X\x Y+X\x Z \ar[r]_-{\psi_{X,Y,Z}} & X\x(Y+Z) }$$
and now combining these we get a commutative diagram
$$\xymatrix @C4pc {
X\x Y+X\x Z \ar@{=}[d] \ar[r]^{\psi_{X,Y,0}+\psi_{X,0,Z}} & 
X\x Y+X\x Z. \ar[d]^{\delta_{X,Y,Z}} \\
X\x Y+X\x Z \ar[r]_-{\psi_{X,Y,Z}} & X\x(Y+Z) }$$
In this last diagram, the $\psi$',s are all invertible, hence so is $\delta$.
\endproof

Recall that an object $T$ is called {\em subterminal} if for any object $X$
there is at most one morphism from $X$ to $T$. (If, as here, a terminal object
exists, this is equivalent to saying that the unique map $T\to 1$ is a
monomorphism.)

Thus to prove our result about non-canonical distributivity
isomorphisms, we must show that the assumption that $X\times 0\cong 0$ 
made in the Lemma is unnecessary. The remainder of this section will
be devoted to doing so.


\begin{proposition}
The product $0\x 0$ is initial, and so $0$ is subterminal.
\end{proposition}

\proof
For the first part, observe that $\psi_{0,0,1}$ gives an isomorphism 
$0\x 0+0\x1\cong0\x(0+1)$, and that $0+1\cong 1$ and $0\x 1\cong 0$.
For the second, we have an isomorphism $0\cong 0\x 0$, and since $0$ 
is initial, this can only be the diagonal $\Delta:0\to 0\x 0$. Thus any morphism
$X\to 0\x 0$ factorizes through the diagonal, and so any two morphisms
$X\to 0$ are equal.
\endproof

Next we consider the special case where \D is pointed $(0=1)$; ultimately we shall 
reduce the general case to this. In a pointed category, every object has a 
(unique) morphism into $0$; but in a distributive category, any morphism into 
$0$ is invertible \cite[Proposition~3.4]{ext}. It follows that any
category which is pointed and distributive is equivalent to the
terminal category 1. Our next result shows that the same conclusion
holds under the assumption of pointedness and a non-canonical
distributivity isomorphism.
 
\begin{proposition}\label{prop:pointed-case}
If \D is pointed  then \D is equivalent to the terminal category 1.
\end{proposition}

\proof
Taking $Y=Z=1$ gives a natural family $\theta_X=\psi_{X,1,1}:X+X\cong X$. By
naturality, the diagram 
$$\xymatrix @C3pc {
X+X+X+X \ar[r]^-{\theta_X+\theta_X}  \ar[d]_{\theta_{X+X}} & X+X \ar[d]^{\theta_X} \\
X+X \ar[r]_-{\theta_X} & X}$$
commutes, and now since $\theta_X$ is invertible $\theta_X+\theta_X=\theta_{X+X}$.
The diagram
$$\xymatrix{
X+X \ar[dd]_{\nabla} \ar[rr]^{\theta_X} \ar[dr]_{i_{X+X}} && X \ar[dr]_{i} \\ 
& X+X+X+X \ar[rr]^{\theta_X+\theta_X}_{\theta_{X+X}} \ar[dd]_{\nabla+\nabla} && 
X+X \ar[dd]_{\nabla} \\
X \ar[dr]_{i} \\
& X+X \ar[rr]_{\theta_X} && X }$$
commutes, where $i_{X+X}$ denotes the injection of the first two copies of $X$
into $X+X+X+X$. But now $\theta_X=\nabla i\theta_X=\theta_X i\nabla$ is invertible,
so $\nabla:X+X\to X$ is a monomorphism; since it also has a section $i$ (and $j$)
it is invertible. This proves that any two maps $X\to Y$ must be equal. On the other hand, there is always at 
least one such map, since \D is pointed; thus there is exactly one, and so
$X\cong 0$. Since $X$ was arbitrary, the result follows.
\endproof

\begin{theorem}
If \D is a category with finite products and coproducts, and with a 
natural family $$\psi_{X,Y,Z}:X\x Y+X\x Z\cong X\x(Y+Z)$$ of isomorphisms,
then \D is distributive.
\end{theorem}

\proof
By Lemma~\ref{lemma:pres+}, it will suffice to show that $X\x 0\cong 0$. Since we have the
projection $X\x 0\to 0$, and the composite $0\to X\x 0\to 0$ is certainly 
the identity, we need only show that the other composite $e:X\x 0\to 0\to X\x 0$
is the identity.  This is an endomorphism in the slice category $\D/0$. So if 
$\D/0$ is trivial, then this composite $e$ will be the identity, and $X\x 0$ will
be isomorphic to $0$. 

Since  $0$ is subterminal, the projection $\D/0\to\D$ is fully faithful, and
preserves finite products as well as coproducts. Thus the isomorphisms $\psi_{X,Y,Z}$
restrict to $\D/0$, thus equipping $\D/0$ with non-standard distributivity 
isomorphisms. By Proposition~\ref{prop:pointed-case}, $\D/0$ is trivial, and so \D is distributive.
\endproof

\section{Non-canonical semi-additivity isomorphisms}
\label{sect:additive}

We now give an analogous result for semi-additivity. An interesting
feature is that this does not require us to assume that the category
is pointed, although that will of course be a consequence.

\begin{theorem}\label{thm:additive}
If \A is a category with finite products and coproducts and with a natural 
family $$\psi_{Y,Z}:Y+Z\cong Y\x Z$$ 
of isomorphisms, then \A is semi-additive. 
\end{theorem}

\proof
Taking $Y=1$ and $Z=0$ gives an isomorphism $\psi_{1,0}:1\cong 1\x 0$;
composing with the projection $1\x 0\to 0$ gives a morphism $1\to 0$. 
By uniqueness of morphisms into $1$ and out of $0$, this is inverse to the 
unique map $0\to 1$, and so \A is pointed. 

Taking one of $Y$ and $Z$ to be $0$ gives natural isomorphisms 
$\psi_{Y,0}:Y\cong Y$ and $\psi_{0,Z}:Z\cong Z$. 
By naturality of the $\psi_{Y,Z}$, the diagrams 
$$\xymatrix @C3pc {
Y \ar[r]^{\psi_{Y,0}} \ar[d]_{i} & Y \ar[d]^{\binom{Y}{0}} &
Z \ar[r]^{\psi_{0,Z}} \ar[d]_{j} & Z \ar[d]^{\binom{0}{Z}} \\
Y+Z \ar[r]_{\psi_{Y,Z}} & Y\x Z & Y+Z \ar[r]_{\psi_{Y,Z}} & Y\x Z }$$
commute, and so also 
$$\xymatrix @C3pc { 
Y+Z \ar[r]^{\psi_{Y,0}+\psi_{0,Z}} \ar@{=}[d] & Y+Z \ar[d]^{\alpha_{Y,Z}} \\
Y+Z \ar[r]_{\psi_{Y,Z}} & Y\x Z }$$
commutes. Just as in the proof of the lemma, $\psi_{Y,0}+\psi_{0,Z}$
and $\psi_{Y,Z}$ are invertible, hence so is $\alpha_{Y,Z}$.
\endproof

\section{Non-canonical isomorphisms for monoidal functors}
\label{sect:monoidal}

In this section we prove a general result on monoidal functors, which
could be used in the proof of both of the other theorems.
Recall that if \A and \B be monoidal categories, a
monoidal
functor $F:\A\to\B$ consists of a functor (also called $F$) equipped
with maps $\phi_{Y,Z}:FY\ox FZ\to F(Y\ox Z)$ and $\phi_0:I\to FI$
which need not be invertible, but which are natural and coherent
\cite{closed}. The monoidal functor is said to be {\em strong}
if $\phi_{Y,Z}$ and $\phi_0$ are invertible, and {\em normal} if 
$\phi_0$ is invertible. Given such an $F$ and another monoidal 
functor $G:\A\to\B$ with structure maps $\psi_{X,Y}$ and $\psi_0$,
a natural transformation $\alpha:F\to G$ is {\em
monoidal} if the diagrams 
$$\xymatrix{
FY\ox FZ \ar[r]^{\alpha_X\ox \alpha_Y} \ar[d]_{\phi_{X,Y}} & GY\ox GZ
\ar[d]^{\psi_{X,Y}} & I \ar[r]^{\phi_0} \ar[dr]_{\psi_0} & FI
\ar[d]^{\alpha I} \\
F(Y\ox Z) \ar[r]_{\alpha_{X\ox Y}} & G(Y\ox Z) && GI }$$
commute. 
Recall further \cite{JoyalStreet-braided} that if \C is braided monoidal, then the functor
$\ox:\C\x\C\to\C$ is strong monoidal, with structure maps 
$$\xymatrix{
W\ox X\ox Y\ox Z \ar[rr]^{W\ox \gamma\ox D} && W\ox Y\ox X\ox Z &
I \ar[r]^-{\lambda} & I\ox I }$$
where $\gamma$ denotes the braiding and $\lambda$ the canonical
isomorphism.

\begin{theorem}\label{thm:monoidal}
Let \A and \B be braided monoidal categories, and $F=(F,\phi,\phi_0):\A\to\B$ a
normal monoidal functor (so that $\phi_0$ is invertible). Suppose further that we have
a monoidal isomorphism
$$\xymatrix @R1pc {
\A\x\A \ar[r]^{F\x F} \ar[dd]_{\ox} & \B\x\B \ar[dd]^{\ox} \\
{}\rtwocell\omit{\psi} & {} \\
\A \ar[r]_{F} & \B }$$
Then $\phi$ is invertible, and so $F$ is strong monoidal.
\end{theorem}

\proof
The fact that $\psi$ is monoidal means in particular that the diagram
$$\xymatrix @C3pc {
FW\ox FX\ox FY\ox FZ \ar[r]^{\psi_{W,X}\ox\psi_{Y,Z}} \ar[d]_{1\ox\gamma\ox 1} &
F(W\ox X)\ox F(Y\ox Z)  \ar[d]^{\phi_{W\ox X,Y\ox Z}} \\
FW\ox FY\ox FX\ox FZ \ar[d]_{\phi_{W,Y}\ox\phi_{X,Z}} & F(W\ox X\ox
Y\ox Z) \ar[d]^{F(1\ox\gamma\ox 1)} \\
F(W\ox Y)\ox F(X\ox Z) \ar[r]_{\psi_{W\ox Y,X\ox Z}} & F(W\ox Y\ox
X\ox Z) }$$
commutes. Taking $X=Y=I$ and twice using the isomorphism $\phi_0$ gives 
commutativity of 
$$\xymatrix{
FW\ox FZ \ar[d]_{1\ox\phi_0\ox\phi_0\ox 1} \ar@{=}[dr] \\
FW\ox FI\ox FI\ox FZ \ar[r]_{\psi_{W,I}\ox\psi_{I,Z}}
\ar[d]_{\phi_{W,I}\ox\phi_{I,Z}} & FW\ox FZ \ar[d]^{\phi_{W,Z}} \\
FW\ox FZ \ar[r]_{\psi_{W,Z}} & F(W\ox Z)}$$
in which all arrows except $\phi_{W,Z}$ are invertible; thus
$\phi_{W,Z}$ too is invertible.
\endproof

\begin{remark}
In the proof of Theorem~\ref{thm:monoidal}, we have used rather less
than was assumed in the statement. For example, we do not use the 
nullary part of the assumption that the natural transformation is
monoidal. 
\end{remark}
 
The following corollary appeared (in dual form) as
\cite[Theorem~3.3]{CaccamoWinskel}:

\begin{corollary}[Caccamo-Winskel]
Let \A and \B be categories with finite coproducts, and $F:\A\to\B$
a functor which preserves the initial object. If there is a natural
family of  isomorphisms
$$\xymatrix{
FX+FY \ar[r]^{\psi_{X,Y}} & F(X+Y) }$$
then $F$ preserves finite coproducts.
\end{corollary}

\proof
In this case $F$ has a unique monoidal structure, and $\psi$ is always
monoidal.
\endproof

In particular if \D has finite products and coproducts, we may apply 
the Corollary to the functor $X\x-:\D\to\D$ and recover 
Lemma~\ref{lemma:pres+}.

Section~\ref{sect:additive} involves the case where the categories \A
and \B are the same, but the monoidal structure on \A is cartesian
and that on \B is cocartesian. The functor $F$ is the identity. One
proves $0\to 1$ is invertible, as in the proof of
Theorem~\ref{thm:additive}; and then the identity $1:\A\to\A$ has a 
unique normal monoidal structure, with binary part precisely the 
canonical morphism $\alpha:Y+Z\to Y\times Z$. Furthermore,
any natural isomorphism $\psi_{Y,Z}:Y+Z\cong Y\x Z$ is monoidal.

\bibliographystyle{plain}

\end{document}